%% file: turan-groups.tex
\documentclass[a4paper]{article}

\usepackage{graphicx}
\usepackage{amsmath}
\usepackage{amsthm}
\usepackage{amssymb}
\usepackage{amsfonts}
\usepackage{latexsym}
\usepackage{epsfig}
\usepackage{color}

\newcommand{\beql}[1]{\begin{equation}\label{#1}}
\newcommand{\eeq}{\end{equation}}
\newcommand{\comment}[1]{}

\newcommand{\eref}[1]{{\rm (\ref{#1})}}

\newcommand{\Abs}[1]{{\left|{#1}\right|}}

\newcommand{\Norm}[1]{{\left\|{#1}\right\|}}

\newcommand{\Qed}{\ \\\mbox{$\Box$}}

\newcommand{\Set}[1]{{\left\{{#1}\right\}}}

\newcommand{\RR}{{\mathbb R}}
\newcommand{\CC}{{\mathbb C}}
\newcommand{\ZZ}{{\mathbb Z}}
\newcommand{\NN}{{\mathbb N}}
\newcommand{\TT}{{\mathbb T}}

\newcommand{\inner}[2]{{\left\langle #1, #2 \right\rangle}}

\newcommand{\dens}{{\rm dens\,}}

\newcommand{\supp}{{\rm supp\,}}

\newcommand{\ft}[1]{\widehat{#1}}

\newcounter{open}
\newcounter{dfn}
\def\thedfn{\arabic{dfn}}

\newcounter{obs}
\def\theobs{\arabic{obs}}

\newcounter{thm}
\newcounter{othm}
\newcounter{mysec}
\def\themysec{\arabic{mysec}}
\newcommand{\mysection}[1]{
  \vskip 0.25in
  \refstepcounter{mysec}\centerline{\large\bf \S\themysec.\ {#1}}\par
  \nopagebreak
  \addcontentsline{toc}{section}{{\bf \themysec.}\ {#1}}
}
\newcounter{mysubsec}[mysec]
\def\themysubsec{\arabic{mysec}.\arabic{mysubsec}}
\newcommand{\mysubsection}[1]{
  \vskip 0.125in
  \refstepcounter{mysubsec}\noindent{\bf \themysubsec\ \ \ {#1}}\par
\nopagebreak
  \addcontentsline{toc}{subsection}{\themysubsec.\ {#1}}
}

\newtheorem{theorem}{Theorem}
\newtheorem{corollary}{Corollary}
\newtheorem{lemma}{Lemma}
\newtheorem{proposition}{Proposition}
\newtheorem{definition}{Definition}

\theoremstyle{definition}
\newtheorem{example}{Example}
\newtheorem{remark}{Remark}

\newcommand{\Tu}{{\mathcal T}}
\newcommand{\Om}{{\Omega}}
\newcommand{\FF}{{\mathcal F}}
\newcommand{\de}{\delta}

\newcounter{rem}
\newcounter{rep}
\newcounter{rev}
\reversemarginpar
\begin{document}

\title{Tur\'an's extremal problem for\\
positive definite functions on groups\footnote{This joint work was done during 
the first author's visit at the Alfr\'ed R\'enyi Institute of Mathematics - 
 Hungarian Academy of Sciences, supported by the ``Mathematics in Information
 Society'' project in the framework of the European Community's ``Confirming
 the International Role of Community Research'' programme.}}
\author{Mihail N. Kolountzakis\footnote{
Supported in part {\it European Commission}
Harmonic Analysis and Related Problems 2002-2006 IHP Network
(Contract Number: HPRN-CT-2001-00273 - HARP)} \and Szil\' ard Gy. R\' ev\' esz
\thanks {The second author was supported in part by the
Hungarian National Foundation for Scientific Research,
Grant \# T034531 and T 032872.}
}
\date{}
\maketitle


\begin{abstract}
We study the following question: Given an open set $\Omega$,
symmetric about $0$, and a continuous, integrable, positive
definite function $f$, supported in $\Omega$ and with $f(0)=1$,
how large can $\int f$ be? This problem has been studied so far
mostly for convex domains $\Omega$ in Euclidean space. In this
paper we study the question in arbitrary locally compact abelian
groups and for more general domains. Our emphasis is on finite
groups as well as Euclidean spaces and $\ZZ^d$. We exhibit upper
bounds for $\int f$ assuming geometric properties of $\Omega$ of
two types: (a) packing properties of $\Omega$ and (b) spectral
properties of $\Omega$. Several examples and applications of the
main theorems are shown. In particular we recover and extend
several known results concerning convex domains in Euclidean
space. Also, we investigate the question of estimating $\int_{\Om}f$ 
over possibly dispersed sets solely in
dependence of the given measure $m:=|\Om|$ of $\Om$. In this respect we 
show that in $\RR$ and $\ZZ$ the integral is maximal for intervals.
\end{abstract}




\mysection{Introduction} \label{sec:introduction}

\mysubsection{A problem of Tur\'an and Stechkin}
\label{sec:a-problem}

\noindent
We study the following problem proposed by Tur\' an and Stechkin \cite{stechkin:periodic}:
\begin{quotation}
Given an open set
$\Omega$, symmetric about $0$, and a continuous,
positive definite, integrable function $f$, with $\supp f \subseteq \Omega$
and with $f(0)=1$, how large can $\int f$ be?
\end{quotation}
The cases studied so far concern $\Omega$ being a convex subset
of $\RR^d$
\cite{arestov:hexagon,arestov:tiles,gorbachev:ball,kolountzakis:turan}
or an interval in the torus $\TT=\RR/\ZZ$
\cite{gorbachev:turan,stechkin:periodic,AKP}.

Such a question is interesting in the study of sphere packing
\cite{gorbachev:sphere,cohn:packings}, in additive number theory
\cite{Ruzsa,KMF} and in the theory of Dirichlet characters and
exponential sums \cite{KS}, among other things.

In this paper we study the problem in more general locally compact
abelian (LCA) groups. This simplifies and unifies many of the
existing results and gives several new estimates and examples. If
$G$ is a LCA group a continuous function $f\in L^1(G)$ is
positive definite if its Fourier transform $\ft{f}:\ft{G}\to\CC$
is everywhere nonnegative on the dual group $\ft{G}$, see
\S\ref{sec:positivedefiniteness}. For the relevant definitions of
the Fourier transform we refer to \cite[Chapter
VII]{katznelson:book} or \cite{rudin:groups}.

The set $\Omega$ will always be taken in this paper to be a
$0$-symmetric, open set in $G$.

If $f\in L^1(G)$ is continuous, positive definite and
supported in $\Omega$ it follows that $f(0) \ge f(x)$ for any
$x\in G$. This leads to the estimate $\int_G f \le \Abs{\Omega}
f(0)$, called the {\em trivial estimate} from now on.

\begin{definition}\label{def:turan-constant}
The {\em Tur\'an constant} ${\cal T}_G(\Omega)$ of a
$0$-symmetric, open subset $\Omega$ of a LCA group $G$ is the
supremum of the quantity $ {\int_G f} /{f(0)} $, where $f\in
L^1(G)$ is continuous and positive definite, and $\supp f$
is a closed set contained in $\Omega$.
\end{definition}

\begin{remark}
The quantity ${\mathcal T}_G(\Omega)$ depends on which
normalization we use for the Haar measure on $G$. If $G$ is
discrete we use the counting measure and if $G$ is compact and
non-discrete we normalize the measure of $G$ to be 1.
\end{remark}

The {\em trivial upper estimate} for the Tur\'an constant is
${\cal T}_G(\Omega) \le \Abs{\Omega}$.

\mysubsection{Previous work}
\label{sec:previous}

\noindent
Let us review some of the known results.

Stechkin \cite{stechkin:periodic} proves $\Tu_{\TT}(\Omega) =
{1\over 2}\Abs{\Omega}$ if $\Omega\subseteq\TT=\RR/\ZZ$ is a
$0$-symmetric interval whose half length divides the length of
$\TT$.

In \cite{arestov:hexagon,arestov:tiles} Arestov and Berdysheva
prove that if $\Omega\subseteq\RR^d$ is a convex polytope which
can tile space when translated by the lattice
$\Lambda\subseteq\RR^d$ (this means that the copies
$\Omega+\lambda$, $\lambda\in\Lambda$, are non-overlapping and
almost every point in space is covered) then $\Tu_{\RR^d}(\Omega)
\le 2^{-d}\Abs{\Omega}$.

Gorbachev \cite{gorbachev:ball} also shows the same inequality
if $\Omega$ is the Euclidean ball
in $\RR^d$ (a different proof of this is given in \cite{kolountzakis:turan}).
The ball clearly cannot tile space.

Kolountzakis and R\'ev\'esz \cite{kolountzakis:turan} show the same inequality for all
convex domains in $\RR^d$ which are spectral (the definition appears later in
this paper in \S\ref{sec:spectral})--convex spectral sets are conjecturally the same
as convex tiles \cite{fuglede:conjecture}.
It is known that all convex tiles are spectral (see e.g.\ \cite{kolountzakis:turan}),
so the result of Arestov and Berdysheva
\cite{arestov:hexagon,arestov:tiles} is also a consequence of the
result in \cite{kolountzakis:turan}.

Gorbachev and Manoshina \cite{gorbachev:turan} study the function
$\Tu_{\TT}(\Omega)$ when $\Omega$ is a $0$-symmetric interval
whose half length does not divide the length of $\TT$, and they
give more detailed information on $\Tu_{\TT}(\Omega)$ when that
length is of a certain arithmetical type.

\mysubsection{Various forms of the Tur\'an problem}

\noindent In fact, it is worth noting that Tur\'an type problems
can be, and have been considered with various settings, although the
relation of these has not been fully clarified yet.
Thus in extending the investigation to LCA groups or to domains in Euclidean
groups which are not convex, the issue of
equivalence has to be dealt with.
One may consider the following function classes.
\begin{eqnarray}
\label{Flclosed} \FF_1(\Om) &:=& \bigg\{f\in L^1(G)~: ~~\supp f
\subset \Om,~~ f \,{\rm positive \,\,
definite}\, \bigg\}\,, \\
\label{Flcontclosed} \FF_{\&}(\Om) &:=& \bigg\{f\in L^1(G)\cap
C(G)~: ~~\supp f \subset \Om,~~ f \,{\rm positive \,\,
definite}\, \bigg\}\,, \\
\label{Flcompact} \FF_c(\Om) &:=& \bigg\{f\in L^1(G)~: ~~\supp f
\subset \subset \Om,~~ f \,{\rm positive \,\,
definite}\, \bigg\}\,, \\
\label{Fcontcompact} \FF(\Om) &:=& \bigg\{f\in C(G)~: ~~\supp f
\subset \subset \Om,~~ f \,{\rm positive \,\, definite}\,
\bigg\}\,\,.
\end{eqnarray}
In $\FF_1, \FF_{\&}$ $\supp f$ is assumed to be merely closed ad
not necessarily compact, and in $\FF_1, \FF_c$ the function $f$
may be discontinuous.

The respective Tur\'an constants are
\begin{eqnarray}\label{Turanconstants}
\Tu_G^{(1)}(\Om)~ {\rm or} ~ \Tu_G^{\&}(\Om)~ {\rm or} ~
\Tu_G^c(\Om)~ {\rm or} ~\Tu_G(\Om) := \qquad \qquad\qquad\qquad
\qquad\qquad\qquad\\
\qquad \qquad\qquad \sup \bigg\{\frac{\int_G f}{f(0)} \,:~ f \in
\FF_1(\Om) ~ {\rm or} ~ \FF_{\&}(\Om) ~ {\rm or} ~ \FF_c(\Om) ~
{\rm or}~ \FF(\Om), ~ {\rm resp.} \bigg\}.\notag
\end{eqnarray}

In general we should consider functions $f:G\to \CC$. But
according to \eqref{posdefdef} also $\overline{f}$ and thus even
$\varphi:=\Re f$ is positive definite, while belonging to the same
function class. As we also have $f(0)=\varphi(0)$ and $\int
f=\int \varphi$, restriction to real valued functions does
not change the values of the Tur\'an constants.

To start with, we prove in \S\ref{sec:positivedefiniteness}

\begin{theorem}\label{th:equivalences}
We have for any LCA group the equivalence of the above defined
versions of the Tur\'an constants:
\begin{equation}\label{equivalance:Turanconstants}
\Tu_G^{(1)}(\Om)=\Tu_G^{\&}(\Om)= \Tu_G^c(\Om)=\Tu_G(\Om)\,.
\end{equation}
\end{theorem}

Note that the original formulation, presented also above in \S
\ref{sec:a-problem} corresponds to $\Tu_G^{\&}(\Om)$.

\begin{remark} It is not fully clarified what happens for
functions vanishing only outside of $\Om$, but having nonzero
values up to the boundary $\partial\Om$.
\end{remark}

\mysubsection{New results} \label{sec:new-results}

\noindent In this paper we focus mostly (but not exclusively) on
finite or compact abelian groups.

Especially in the case of finite groups we can show clearly the
geometric aspects of the problem without being sidetracked by
technicalities that arise when the group is not discrete or not
compact.

We present two types of results. In the first type some kind of
``packing'' condition is assumed on $\Omega$ which leads to an
upper bound for ${\cal T}_G(\Omega)$. (The justification of the
term ``packing'' should be more evident in the statement of
Corollary \ref{cor:packing} in \S\ref{sec:packing-proofs}.)
\begin{theorem}\label{th:diff}
Suppose that $G$ is a compact abelian group, $\Lambda\subseteq G$,
$\Omega\subseteq G$ is a $0$-symmetric open set and
$(\Lambda-\Lambda) \cap \Omega \subseteq \Set{0}$. Suppose also
that $f\in L^1(G)$ is a continuous positive definite function
supported on $\Omega$. Then
\begin{equation}\label{bound-from-diff}
\int_G f(x)\,dx \le {\Abs{G} \over \Abs{\Lambda}} f(0).
\end{equation}
In other words ${\cal T}_G(\Omega) \le \Abs{G} / \Abs{\Lambda}$.
\end{theorem}
(Observe that the conditions imply that $\Lambda$ is finite.)

The proof appears in \S\ref{sec:upper-bound-from-packing}.

The following Theorem \ref{th:diff-infinite} is analogous to
Theorem \ref{th:diff} for the non-compact case.

\begin{theorem}\label{th:diff-infinite}
Suppose that $G$ is one of the groups $\RR^d$ or $\ZZ^d$, that
$\Lambda \subseteq G$ is a set of upper density $\rho>0$, and
$\Omega\subseteq G$ is a $0$-symmetric open set such that $\Omega
\cap (\Lambda-\Lambda) \subseteq \Set{0}$. Let also $f\in L^1(G)$
be a continuous positive definite function on $G$ whose support
is a compact set contained in $\Omega$. Then
\begin{equation}\label{bound-from-diff-infinite}
\int_G f(x)\,dx \le {1 \over \rho} f(0).
\end{equation}
In other words ${\cal T}_G(\Omega) \le 1/\rho$.
\end{theorem}

In \S\ref{sec:sub-groups} and
\S\ref{sec:sharpness}--\S\ref{sec:dispersed} we present several
examples and applications of Theorems \ref{th:diff} and
\ref{th:diff-infinite} in various groups. These theorems in
particular imply the results of
\cite{stechkin:periodic,arestov:hexagon,arestov:tiles}, but are
much more general.

The second type of result we give is analogous
to that proved in \cite{kolountzakis:turan}.
Here we suppose that $\Omega$ can be embedded in the difference set
of a spectral set (see definition in \S\ref{sec:spectral})
and we derive an upper bound for ${\cal T}_G(\Omega)$ from that.
\begin{theorem}\label{th:spectral}
Suppose $G$ is a finite abelian group, $\Omega,H \subseteq G$,
$\Omega \subseteq H-H$, and that $H$ is a spectral set with
spectrum $T\subseteq\ft{G}$. Then for any positive definite
function on $G$ with support in $\Omega$ we have
\begin{equation}\label{bound-from-spectrum}
\sum_{x\in G}f(x) \le \Abs{H} f(0).
\end{equation}
In other words ${\cal T}_G(\Omega) \le \Abs{H}$.
\end{theorem}

What was essentially proved in \cite{kolountzakis:turan} was a
``continuous'' version of Theorem \ref{th:spectral}. Essentially,
the following was proved.

\begin{theorem}\label{th:spectral-infinite}{\rm \cite{kolountzakis:turan}}
If $H$ is a bounded open set in $\RR^d$ which is spectral, then
for the difference set $\Om=H-H$ we have ${\cal T}_{\RR^d}(\Om) =
\Abs{H}$.
\end{theorem}

The result there was only formulated for convex sets
$H\subset\RR^d$ (for which $H-H=2H$) but the proof works verbatim
for the result we just stated. Let us emphasize here that in the
continuous case we demand that eligible functions for our
extremal problem have compact support contained in the given open
set $\Omega$ (whose Tur\'an constant we are estimating).

We give the proof of Theorem \ref{th:spectral} in
\S\ref{sec:upper-bound-from-spectral-sets}.

Furthermore, in \S\ref{sec:comparison} we show that there are
cases when Theorems \ref{th:spectral} and
\ref{th:spectral-infinite} give provably better results than any
application of Theorems \ref{th:diff} and \ref{th:diff-infinite},
respectively. For this we use one of Tao's
\cite{tao:fuglede-down} recent examples which show one direction
of Fuglede's conjecture to be false.


\mysection{Preliminaries}
\label{sec:preliminaries}

\noindent In this section we describe the basic facts about
positive definite functions, translational tiling, packing and
spectral sets on LCA groups.


\mysubsection{Positive definite functions on LCA groups}
\label{sec:positivedefiniteness}

\noindent In this subsection we explore a few facts on positive
definite, not necessarily continuous functions.

Recall that on a LCA group $G$ a function $f$ is called positive
definite if the inequality
\begin{equation}\label{posdefdef}
\sum_{n,m=1}^{N}~ c_n \overline{c_m} f(x_n-x_m)\ge 0
\qquad (\forall x_1,\dots,x_N\in G, \forall c_1,\dots,c_N\in\CC)
\end{equation}
holds true. Note that positive definite functions are not
assumed to be continuous. Still, all such functions $f$ are
necessarily bounded by $f(0)$ \cite[p.\ 18, Eqn (3)]{rudin:groups}.
Moreover,
$f(-x)=\widetilde{f}(x):=\overline{f(x)}$ for all $x\in G$
\cite[p.\ 18, Eqn (2)]{rudin:groups}, hence the support of $f$ is
necessarily symmetric, and the condition $\supp f \subset \Om$
implies also $\supp f \subset \Om\cap(-\Om)$. The latter set
being symmetric, without loss of generality we can assume at the
outset that $\Om$ is symmetric itself.

It is immediate from \eqref{posdefdef} that for any subgroup $K$
of $G$, the restriction $f|_K$ of a positive definite function
$f$ is also positive definite on $K$.

The Fourier transform $\ft{f}$ of an $f\in L^1(G)$ belongs to
$A(\ft{G})\subset C_0(\ft{G})$, and the Fourier transform of the
convolution $f*g$ of $f,g\in L^1(G)$, defined almost everywhere,
satisfies $\ft{f*g}=\ft{f}\ft{g}$ \cite[Theorem 1.2.4]{rudin:groups}.
Similarly, for $\nu,\mu\in M(G)$ and their convolution
$\mu*\nu\in M(G)$ the Fourier transforms are bounded and uniformly
continuous and $\ft{\mu*\nu}=\ft{\mu}\ft{\nu}$ \cite[Theorem 1.3.3]{rudin:groups}.

In case $f,g\in L^2(G)$, the convolution $h:=f*g$ is defined even
in the pointwise sense and $h\in C_0(\widehat{G})$ \cite[Theorem 1.1.6(d)]{rudin:groups}.
For $f\in L^2(G)$ arbitrary (denoting as above,
$\widetilde{f}(x):=\overline{f(-x)}$), $f*\widetilde{f}$ is
continuous and positive definite with Fourier transform
$|\ft{f}|^2$ \cite[1.4.2(a)]{rudin:groups}.

Note that for any given $\gamma\in\ft{G}$ $f$ is positive definite
if and only if $f(x)\gamma(x)$ is positive definite; this can be
checked by modifying the coefficients in \eqref{posdefdef}
accordingly.

\begin{lemma}\label{lemma:productposdef} Suppose that $f$
is (measurable and) positive definite and $g\in L^2(G)$ is
arbitrary. Then the product $f\cdot(g*\widetilde{g})$ is positive
definite.
\end{lemma}

\begin{proof} As written above, $h:=g*\widetilde{g}\in C_0(G)$, while
$f$, being positive definite, is also bounded. Take now $x_n\in
G$ and $c_n\in\CC$ for $n=1,\dots,N$ arbitrarily. Then
\begin{eqnarray*}
& & \sum_{n,m=1}^{N}~ c_n \overline{c_m} f(x_n-x_m)h(x_n-x_m) \\
&=& \sum_{n,m=1}^{N}~ c_n \overline{c_m} f(x_n-x_m) \int_G
g(x_n-y)\overline{g}(x_m-y)dy \\
&=& \int_G \sum_{n,m=1}^{N}~ a_n(y) \overline{a_m(y)} f(x_n-x_m)
~dy\,,
\end{eqnarray*}
where $a_n(y):=c_n g(x_n-y)\in L^2(G)$ ($n=1,\dots,N$). Since the
expression under the integral sign is nonnegative by
\eqref{posdefdef} for each given $y$, also the integral is
nonnegative and the assertion follows.
\end{proof}

Note that we did not assume $f$ to be integrable, and neither the
product $fh$ is supposed to belong to any subspace. By positive
definiteness, $f$ is bounded; but if $G$ is not compact, $\ft{f}$
is not necessarily defined. However, as $h\in C_0(G)$, in any case
we must have $fh \in L^{\infty}(G)$. This follows from positive
definiteness of $fh$, too.

The next Lemma is obvious for compact groups as we can take
$k=1$.
\begin{lemma}\label{lemma:posdefmultiplier}
Suppose $C$ is a compact set in a LCA group $G$ and $\delta>0$ is
given. Then there exists a compactly supported, positive definite
and continuous ``kernel function'' $k(x)\in C_c(G)$ satisfying
$k(0)=1$, $0\le k \le 1$, and $k|_C\ge 1-\delta$. Moreover, we
can take $k=h*\widetilde{h}$, where $h$ is the $L^2$-normalized
indicator function of a suitable Borel measurable set $V\supset
C$ with compact closure $\overline{V}$.
\end{lemma}

\begin{proof}
We may clearly assume that $G$ is not compact.

The deduction will follow the proof of 2.6.7 Theorem on page 52
of \cite{rudin:groups} with a slight modification towards the end
of the argument. In this proof the compact set $C$ is given, and
then another Borel set $E$ and an increasing sequence of Borel
sets $V_N \quad (N\in \NN)$ are found, so that $C\subset E=V_0$
and $|V_N|=(2N+1)^{n}|E|$
(with $n$ a fixed nonnegative integer constant); moreover, all
the $V_N$ have compact closure and $V_N+E\subset V_{N+s}$ is
ensured for some fixed $s$ and for all $N\in\NN$. Hence for every
$c\in C\subset E$ we have $V_{N+s}-c\supset V_N$. Denoting the
indicator function of $V_{N+s}$ by $\chi$ we are led to $\int_G
\chi(x+c)\chi(x)dx\ge |V_N|$. Putting $h:=|V_{N+s}|^{-1/2}\chi$
yields $h*\widetilde{h}(c)\ge |V_N|/|V_{N+s}|>1-\delta$, if $N$ is
chosen large enough (depending on the constants $n$, $s$ and the
given $\delta$). With this choice of $h$ and $V:=V_{N+s}$ all
assertions of the Lemma are true.
\end{proof}

\begin{remark} As Rudin points out, this argument essentially depends on
structure theorems of LCA groups.
\end{remark}

\begin{lemma}\label{lemma:posdef-ftnonnegative} Suppose that
$f\in L^1(G)$ is positive definite. Then the Fourier transform
$\ft{f}$ is nonnegative.
\end{lemma}

\begin{proof} Since
for any character $\gamma\in\ft{G}$ we have $\ft{\gamma f}
=\ft{f}(\cdot-\gamma)$, and $f$ is positive definite precisely
when $\gamma f$ is such, it suffices to prove that $\ft{f}(0)\ge
0$.

For technical reasons, we need to modify $f$ to have compact
support. Let $\de$ be any positive parameter. Since
$d\nu(x):=f(x)dx$ is a {\it regular} Borel measure, for some
compact set $C$ we have $\|f\|_{L^(G\setminus C)}<\de$. Take the
function $k$ provided by Lemma \ref{lemma:posdefmultiplier} for
the compact set $C$ and the chosen parameter $\de>0$. If $g:=kf$,
Lemma \ref{lemma:productposdef} shows that $g$ is positive
definite, while $|\ft{g}(0)-\ft{f}(0)| \le |\int_C f - \int_C g| +
\|f\|_{L^1(G\setminus C)} < \de \int_C |f| + \de \le
\de(1+\|f\|_{L^1})$. Choosing $\de$ small enough, it follows that
there exists a compactly supported positive definite $g\in
L^1(G)$ with $\ft{g}(0)<0$ provided that $\ft{f}(0)<0$. Hence it
suffices to prove the assertion for compactly supported positive
definite functions $g$.

Applying definition \eqref{posdefdef} with
all $c_n$ chosen as 1 yields
$$
0\le\sum_{n=1}^{N}\sum_{m=1}^{N}~ g(x_n-x_m)\,.
$$
Integrating over $C^N$ (where $C:=\supp g$) we obtain
$$
0\le N |C|^N g(0) + (N^2-N)|C|^{N-1}\int_C g\,,
$$
which implies
$$
-\frac{|C|g(0)}{N-1}\le \ft{g}(0).
$$
Letting $N\to\infty$ concludes the proof.
\end{proof}

\begin{lemma}\label{lemma:convolutionposdef} Suppose that
$f,g \in L^1(G)$ are two positive definite functions. Then the
convolution $f*g\in L^1(G)$ is uniformly continuous and positive
definite.
\end{lemma}

\begin{proof} Since a positive definite function is bounded, we
have also $f\in L^{\infty}(G)$, hence $f*g$ is uniformly
continuous c.f. \cite[Theorem 1.1.6(b)]{rudin:groups}. For the
Fourier transform $\ft{f*g}=\ft{f}\ft{g}$ of the continuous
function $f*g$ positive definiteness is equivalent to
$\ft{f*g}\ge 0$. Now Lemma \ref{lemma:posdef-ftnonnegative} gives
$\ft{f}\ge0$ and $\ft{g}\ge 0$, hence $\ft{f*g}\ge 0$ and $f*g$
is positive definite.
\end{proof}

\begin{lemma}\label{lemma:posdefunit}
Suppose $U$ is a given neighborhood of $0$ in a LCA group $G$.
Then there exists a compactly supported, continuous, positive
definite and nonnegative ``kernel function'' $k(x)\in C_c(G)$
satisfying $\supp k \subset \subset U$ and $\int k =1$. Moreover,
we can take $k=h*\widetilde{h}$, with $h=|W|^{-1}\chi_W$, where
$\chi_W$ is the indicator function of a compact set $W$ satisfying
$W-W \subset U$.
\end{lemma}

\begin{proof}
By continuity of the operation of subtraction, there exists a
compact neighborhood $W$ of $0$ satisfying $W-W\subset U$. With
the above definitions of $k$ and $h$ we clearly have $\supp k
\subset \subset W-W \subset U$ (c.f. \cite[Theorem
1.1.6(c)]{rudin:groups} and also $\int k = |W|^{-2}(\int
\chi_W)^2 =1$. Since $h, \widetilde{h} \in L^2(G)$, $k\in
C_0(G)$, and as $\supp k$ is compact, $k \in C_c(G)$. Since $h$
is nonnegative, so is $f$. Finally, \cite[1.4.2(a)]{rudin:groups}
gives positive definiteness of $k$.
\end{proof}

\begin{lemma}\label{lemma:convolution-approximation}
Let $f$ be positive definite and integrable. Then for any
$\epsilon>0$ and open set $U$ containing $0$, there exists a
nonnegative, positive definite function of the form
$k=h*\widetilde{h}$ (with $h\in L^2(G)$), so that $\supp k
\subset \subset U$, $\int_U k=1$, and $\|f-f*k\|_1<\epsilon$.
\end{lemma}

\begin{proof} For the given function $f$ there exists a
neighborhood $V$ of $0$ with the property that
$\|f-f*u\|<\epsilon$ whenever $\int_G u=1$ and $u\ge 0$ is Borel
measurable and vanishing outside $V$ \cite[Theorem 1.1.8]{rudin:groups}.
Now we can construct for the open set $U_0:=V \cap
U$ the kernel function $k$ as in Lemma \ref{lemma:posdefunit}.
Clearly, $k$ satisfies all conditions for $u$, hence
$\|f-f*k\|_1<\epsilon$ follows. By construction, $\supp k \subset
\subset U_0 \subset U$ and $\int_U k=1$.
\end{proof}

\begin{lemma}\label{lemma:compact-plus-open}
For any pair of sets $K\subset\subset U$ with $K$ compact and $U$
open, there exists a neighborhood $V$ of $0$ satisfying
$K+V\subset U$.
\end{lemma}

\begin{proof}
Since addition is continuous, for any open neighborhood
$U_0$ of $0$ there exists a neighborhood $W$ so that $W+W\subset
U_0$. Take now to each point $x\in K$ an open neighborhood $W_x$
of $0$ such that $x+W_x+W_x \subset U$, ie. $W_x+W_x\subset U-x$.
Clearly the family of open sets $\{x+W_x~:~x\in K\}$ form an open
covering of $K$, so in view of compactness of $K$ there exists a
finite subcovering $\{W_{x_k}+x_k~:~k=1,\dots,n\}$. Take now
$V:=\bigcap_{k=1}^{n}W_{x_k}$. We claim that $K+V\subset U$.
Indeed, if $y\in K$ and $z\in V$ then considering any index $k$
with $y\in x_k+W_{x_k}$, we find $ y+z \in (x_k+W_{x_k})+V
\subset x_k +W_{x_k} +W_{x_k} \subset U$.
\end{proof}

\begin{lemma}\label{lemma:continuous-approx}
Let $\epsilon>0$ be arbitrary. Assume that $f$ is measurable and
positive definite and compactly supported in the open set $\Om$.
Then there exists another positive definite, but also continuous
function $g$ with $f(0)\ge g(0)$ and $\int_G g\ge \int_G
f-\epsilon $, also supported compactly in $\Om$.
\end{lemma}
\begin{proof}
Observe that $f$, being positive definite, is also bounded, and
since it is compactly supported, it also belongs to $L^1(G)$.
Thus we can use the Fourier transform $\ft{f}$. Let $K:=\supp f
\subset \subset \Om$ and consider a neighborhood $U$ of $0$ with
$K+U\subset\Om$. Such a $U$ is provided by Lemma
\ref{lemma:compact-plus-open}. Lemma
\ref{lemma:convolution-approximation} provides a positive
definite, continuous kernel $k\in C_c(G)$, compactly supported in
$U$ and satisfying $\int_G g\ge\int_G f - \epsilon$. In view of
$k=h*\widetilde{h}$ and Lemma \ref{lemma:convolutionposdef} also
$g:=f*k$ is positive definite while obviously $g\in C_c(G)$ is
supported compactly in $K+U\subset\Om$. It remains to note that
by $k\ge 0$, $\int k=1$ and $\Abs{f}\le f(0)$ we also have $g(0)=\int
k(x)f(-x)dx\le f(0)\int k =f(0)$.
\end{proof}

\begin{proposition}\label{prop:L1} With the definitions above we have
$\Tu_G^{(1)}(\Om)=\Tu_G^c(\Om)$.
\end{proposition}

\begin{proof} Let $\epsilon>0$ and $\de>0$ be arbitrary and
$f\in \FF_{1}(\Om)$ be chosen so that $\int_G
f>\Tu_G^{(1)}(\Om)-\de $. As $f\in L^1(G)$, the measure $|f(x)|dx$
is absolutely continuous with respect to the Haar measure, hence
it is also a regular Borel measure and there exists a compact
subset $C\subset \subset \supp f$ so that $\int_{G\setminus C} |f|
< \de$. Now an application of Lemma \ref{lemma:posdefmultiplier}
with $C$ and $\de$ provides us the positive definite, compactly
supported kernel function $k$ satisfying $k(0)=1$, and $k|_C
>(1-\de)$. Let $g:=fk$. Then $\supp g \subset (\supp k \cap \supp
f) \subset \subset \supp f$, hence $g$ is compactly supported
within $\Om$. Moreover, $g(0)=1$ and $g$ is positive definite in
view of Lemma \ref{lemma:productposdef}. Hence $g\in \FF_c(\Om)$.
We now have
\begin{eqnarray*}
\int g = \int_\Omega kf &=&  \int_\Omega f -  \int_\Omega (1-k)f \\
 &\ge& \int_\Omega  f - \delta\int_C\Abs{f} - \int_{\Omega\setminus C}\Abs{f}\\
 &\ge& \int_\Omega f - \delta \int_{\Om}\Abs{f}- \de
 \qquad \ge (1-\de)\big(\Tu_G^{(1)}(\Om)-\de\big)-\de .
\end{eqnarray*}
Clearly, if $\de$ was chosen small enough, we obtain $\int g >
\Tu_G^{(1)}(\Om)-\epsilon $. Now taking sup over $g \in
\FF_c(\Om)$ concludes the proof, since $\epsilon>0$ was arbitrary.
\end{proof}

\begin{proposition}\label{prop:continuous}
With the definitions above we have $\Tu_G(\Om)=\Tu_G^c(\Om)$.
\end{proposition}
\begin{proof} Since $\FF_c(\Om)\supset\FF(\Om)$, it suffices to
prove $\Tu_G^c(\Om)\le \Tu_G(\Om)$.

Let $\epsilon>0$ and $f\in \FF_c(\Om)$ be chosen so that $\int f
> \Tu_G^c(\Om)- \epsilon$, while $\supp f$ is a compact subset of
the open set $\Om$. Hence an application of Lemma
\ref{lemma:continuous-approx} provides a $g\in\FF(\Om)$ with
$\Tu_G(\Om) \ge \int g > \int f - \epsilon > \Tu_G^c(\Om)-
2\epsilon$. Now $\epsilon\to 0$ yields the Proposition.
\end{proof}

\begin{proof}[Proof of Theorem \ref{th:equivalences}]
We have the obvious inclusions $\FF_{1}(\Om)\supset \FF_{\&}(\Om)
\supset \FF(\Om)$ and $\FF_1(\Om) \supset \FF_c(\Om) \supset
\FF(\Om)$, hence $\Tu_G^{1}(\Om) \ge \Tu_G^{\&}(\Om) \ge
\Tu_G(\Om)$ and $\Tu_G^{(1)}(\Om)\ge \Tu_G^{c}(\Om)\ge
\Tu_G(\Om)$. On combining these inequalities with Propositions
\ref{prop:L1} and \ref{prop:continuous} the assertion follows.
\end{proof}

If we consider a {\it continuous} positive definite function $f$,
then it must also be uniformly continuous \cite[p.\ 18, Eqns (3), (4)]{rudin:groups}.
When $\supp f$ has bounded Haar measure
(and, in particular, when $\supp f$ is compact) then $f$ belongs
to $L^1(G)$, too. For an integrable, continuous and positive
definite function $f$ the Fourier transform $\ft{f}$ of $f$
exists, and the Fourier inversion formula holds, cf.
\cite[1.5.1]{rudin:groups}. The well-known Bochner-Weil
characterization says that $f\in C(G)$ being positive definite is
equivalent to the existence of a non-negative measure $\mu$ on
the dual group $\ft{G}$ so that
$$
f(x)=\int_{\ft{G}} \overline{\gamma(x)}~d\mu(\gamma);
$$
moreover, this representation is unique cf. \cite[1.4.3]{rudin:groups},
Comparing the Fourier inversion formula and the unique
representation above leads to the further characterization that
{\em for a continuous and integrable} $f$ being positive definite
is equivalent to $\ft{f}\ge 0$, compare
\cite[1.7.3(e)]{rudin:groups}. Thus it is really advantageous to restrict
the function class considered from $\FF_{1}(\Om)$ to $\FF(\Om)$,
say.

Our setting is that $\Om$ is an open (symmetric) set, and we
require that $f$ can be nonzero only in $\Om$. This is an
essential condition. In this respect approximation has its
limitations: eg. we cannot relax the conditions to require $\supp
f \subset \overline{\Om}$ only.

Indeed, if $\Om$ is not {\it fat}, meaning that $\Om={\rm int}~
\overline{\Om}$, this can lead to essential changes of the
Tur\'an constants. Eg. if $G=\RR$ and $\Om=(-a,a)\setminus\{\pm
b\}$, then ${\rm int}\overline{\Om}=(-a,a)$ and
$\Tu_{\RR}((-a,a))=a$, while $\Tu_{\RR}(\Om)=b$ if $a/2\le b \le
a$, see Theorem \ref{th:interval-minus-point} below. Similarly,
if $G=\TT$ and $\Om=\TT\setminus \{\pi\}$, then
$\Tu_{\TT}(\Om)=1/2$, but obviously
$\Tu_{\TT}(\overline{\Om})=1$. That is, forcing the function $f$
to vanish at one single point can, through positive definiteness,
bring down the values essentially in general.

In this respect, original formulations of the Tur\'an problem in
\cite{arestov:tiles} and \cite{kolountzakis:turan} may be
misleading, since for a convex body $\Om$ in $\RR^d$ or $\TT^d$
the allegedly extremal function $\chi_{\Om/2}*\chi_{\Om/2}$ does
{\it not} belong to the function class $\FF_{\&}(\Om)$ considered
there. Instead, a corresponding limiting argument should provide
the same extremal value. In convex or star bodies in Euclidean
spaces one can easily obtain a positive definite function
supported properly in the body from one that may be ``non-zero up
to the boundary'', by a slight dilation of space, without losing
much integral. It is unclear how to do this in general, even for
domains in $\RR^d$.

\mysubsection{Tiling and packing}
\label{sec:tiling}

\noindent Suppose $G$ is a LCA group. We say that a function
$f\in L^1(G)$ tiles $G$ by translation with a set
$\Lambda\subseteq G$ at level $c\in\CC$ if
$$
\sum_{\lambda\in\Lambda} f(x-\lambda) = c
$$
for a.a. $x\in G$, with the sum converging absolutely. We then
write ``$f+\Lambda = c G$''.

We say that $f$ {\em packs} $G$ with the translation
set $\Lambda$ at level $c\in\RR$, and write $f+\Lambda \le c G$,
if
$$
\sum_{\lambda\in\Lambda}f(x-\lambda) \le c,
$$
for a.a. $x\in G$.

When the group is finite (and we do not, therefore, have to worry
about the set $\Lambda$ being finite or not) the tiling condition
$f+\Lambda = c G$ means precisely $f*\chi_\Lambda = c$. Taking
Fourier transform, this is the same as $\ft{f}\ft{\chi_\Lambda} =
c\Abs{G}\chi_\Set{0}$, which is in turn equivalent to the
condition
\begin{equation}\label{spectral-condition-for-tiling}
\supp\ft{\chi_\Lambda} \subseteq \Set{0} \cup \Set{\ft{f}=0}
\ \ \mbox{and}\ \ c={\Abs{\Lambda} \over \Abs{G}} \sum_{x\in G} f(x).
\end{equation}

Finally,
if $E \subseteq G$ we say that $E$ packs with $\Lambda$
if $\chi_E$ packs with $\Lambda$ at level $1$.
Observe that $E$ packs with $\Lambda$ if and only if
$$
(E-E) \cap (\Lambda-\Lambda) = \Set{0}.
$$

\mysubsection{Spectra}
\label{sec:spectral}

\noindent Let $G$ be a LCA group and $\ft{G}$ be its dual group,
that is the group of all continuous group homomorphisms
$G\to\CC$. A set $H \subseteq G$ has the set $T \subseteq \ft{G}$
as a {\em spectrum} if and only if $T$ forms an orthogonal basis
for $L^2(H)$.

Suppose from now on that $G$ is finite.

It follows that $\Abs{T} = \Abs{H}$, the dimension of $\ell^2(H)$,
and with a little more work it follows that $T$ is a spectrum of
$H$ if and only if we have the tiling condition
\begin{equation}\label{tiling-condition-for-spectrum}
\Abs{\ft{\chi_H}}^2 + T = \Abs{H}^2 \ft{G}.
\end{equation}
Indeed, for $t_1, t_2 \in \ft{G}$ we have by definition of the Fourier transform
that
$$
\inner{t_1}{t_2}_H = \sum_{x\in H} t_1(x) \overline{t_2}(x) =
 \sum_{x \in H}(t_1-t_2)(x) = \ft{\chi_H}(t_1-t_2).
$$
Suppose now that $T$ is a spectrum of $H$.
If $t\in \ft{G}$ we have (Parseval)
\begin{eqnarray*}
\Abs{H} &=& \Norm{t}_{\ell^2(H)}^2\\
 &=& \sum_{s\in T} \Abs{\inner{t}{{s\over\Norm{s}}}}^2\\
 &=& {1\over\Abs{H}} \sum_{s\in T} \Abs{\inner{t}{s}}^2\\
 &=& {1\over\Abs{H}} \sum_{s\in T} \Abs{\ft{\chi_H}(t-s)}^2,
\end{eqnarray*}
which is precisely the statement that $\Abs{\ft{\chi_H}}^2+T = \Abs{H}^2\ft{G}$.
That this tiling condition is also sufficient to imply that $T$ is a spectrum of
$H$ follows similarly (we are not using this direction in this paper).

By the analysis of tiling shown in \S\ref{sec:tiling} it follows that this happens
if and only if
\begin{equation}\label{spectral-condition-for-spectrum}
\supp\ft{\chi_T} \subseteq \Set{0} \cup (H-H)^c
\ \ \mbox{and}\ \
\Abs{T} = \Abs{H}.
\end{equation}

\mysection{Generalities about Tur\'an constants on groups}
\label{sec:generalities-on-groups}

\mysubsection{Homomorphic images and the Tur\'an constant}
\label{sec:homomorphic-image}

Let $G$ and $H$ be two LCA groups, and $\varphi:G\to H$ a
continuous group homomorphism which maps $G$ onto $H$. Denote
$K:={\rm Ker} (\varphi)\le G $. By continuity of $\varphi$, $K$ is
a closed subgroup, hence a LCA group itself. We consider $G/K$ as
fixed together with the canonical or natural projection $\pi:G\to
G/K$ defined as $\pi(g):=[g]:=g+K\in G/K$. By definition of the
topology of $G/K$, $\pi$ is an open and continuous mapping.
Compare \S B.2, B.6 in Appendix B of \cite{rudin:groups}.
Moreover, $\varphi\circ\pi^{-1}: G/K \to H$ is an isomorphism of
the LCA groups $G/K$ and $H$.

For the determination of the Tur\'an constants, the choice of
the Haar measure is relevant. Haar measures are unique up to
a constant factor: we can always choose the Haar measures $\mu_K$
and $\mu_{G/K}$ so that $d\mu_G=d\mu_K d\mu_{G/K}$, in the sense
of (2) in \cite[2.7.3]{rudin:groups}. On the other hand
fixing a particular Haar measure $\mu_H$ of $H$ always leaves
open the question of compatibility with the fixed measure
$\mu_{G/K}$ and the mapping $\varphi$. Let $A\subset H$ be an
arbitrary Borel set. Then one can define $\nu(A):=\mu_{G/K}
(\pi(\varphi^{-1}(A)))$; since $\varphi$ is onto, clearly this
defines another Haar measure on $H$. Since Haar measures are
constant multiples of each other, we necessarily have
$C:=d\mu_H/d\nu$ a constant. Once $H$ and $\mu_H$ are given,
various homomorphisms $\varphi$ may generate different measures,
but the constant $C=C(\varphi)$ can always be read from this
relation.

\begin{proposition}\label{prop:homomorphism}
Let $G$ and $H$ be LCA groups, and $\varphi : G\to H$ be a
continuous group homomorphism onto $H$. Suppose an open subset
$\Omega\subset G$ is given, and let
$\Theta:=\varphi(\Omega)\subset H$. Consider the closed subgroup
$K:={\rm Ker} (\varphi)\le G $, and the quotient group $G/K$
together with their Haar measures $\mu_{G/K}$ and $\mu_K$,
normalized as above. We then have
\begin{equation}\label{homomorphism}
{\mathcal T}_G(\Omega) \le \frac 1C {\mathcal
T}_H(\Theta){\mathcal T}_K(\Omega\cap K) \qquad (C:= \frac{d
\mu_H}{d\nu })~~.
\end{equation}
Here $\nu:=\mu_{G/K} \circ \pi \circ \varphi^{-1}$ is defined as
above.
\end{proposition}

\begin{proof}
As $K$ is the kernel of the continuous homomorphism $\varphi$, $K$
is a closed subgroup of $G$. Therefore, the factor group $G/K$ is
a LCA group, which is continuously isomorphic to $H$.

The image $\Theta$ of the open set $\Om$ is open, since $\varphi$
is also an open mapping. Indeed, $\pi$ is open by its definition,
and thus $\pi(\Om)$ is open in $G/K$ for any open $\Om$ in $G$.
However, the isomorphism $\psi:G/K\to H$, defined by
$\psi:=\varphi \circ \pi^{-1}$, brings over the open set
$\pi(\Om)$ to $\Theta$, which is then open by the isomorphism
itself.

Observe that $\Om_g:=\Om \cap (K+g)$ is relatively open for any
$g\in G$, while the coset $K+g$ is closed. Let us choose
arbitrarily a representative $g(h)\in G$ of each coset
$\varphi^{-1}(h)$ of $K$ to all $h\in H$. Now for any uniformly
continuous function $f:G\to \CC$ we can define
\begin{equation}\label{projection-integral}
F(h):=\int_K f(g(h)+k)~d\mu_K(k)=\int_{\varphi^{-1}(h)} f(x)~
d\mu_K(x-g(h)).
\end{equation}
Since $f$ is uniformly continuous, the function $F:H\to\CC$ is
continuous, $F(0)=\int_K f d\mu_K$, and by Fubini's Theorem
\begin{align}\label{quotientintegral}
\int_H F(h) d\mu_H(h)& = \int_H\int_{\varphi^{-1}(h)} f(g(h)+k)~
d\mu_K(k)C d\nu(h)
\notag \\
& = C \int_H \int_K f(g(h)+k) ~ d\mu_K(k) d\mu_{G/K}(\pi
\varphi^{-1}(h))
\\
& = C \int_{H\times K} f(g(h)+k) ~ d\mu_K(k) d\mu_{G/K}([g(h)])  =
C \int_G f ~ d\mu_G \,,\notag
\end{align}
taking into account the choice of normalization of the Haar
measures for $K$ and $G/K$. Next we prove that $F$ is positive
definite on $H$ in case $f$ is positive definite on $G$. Indeed,
for any character $\chi$ on $H$ there is a character
$\gamma:=\chi\circ\varphi$ on $G$, and applying
\eqref{quotientintegral} to $f\gamma$ yields
$$
\int_H F(h) \chi(h) d\mu_H(h) = C \int_G f(g)\gamma(g) \ge 0~.
$$
Thus we have $\int_H F d\mu_H\le \Tu_H(\Theta) F(0)$. Moreover,
$f|_K$ is positive definite on $K$, hence we also have
$F(0)=\int_{K\cap\Om}f~d\mu_K \le \Tu_K(K\cap\Om) f(0)$.
Comparing these inequalities with \eqref{quotientintegral} yields
$C \int_G f ~ d\mu_G \le \Tu_H(\Theta) \Tu_K(K\cap\Om) f(0)$, and
taking supremum of $\int_G f d\mu_G/f(0)$ \eqref{homomorphism}
obtains.

\end{proof}

\mysubsection{Automorphic invariance of the Tur\'an constant}
\label{sec:automorphic-images}

One of the reasons to work out Proposition
\ref{prop:homomorphism} is its corollary to the case when we deal
with an automorphism of the group $G$.

\begin{corollary}\label{cor:automorphinvariance}
Let $G$ be a LCA group and let $\varphi:G\to G$ be an automorphism.
Then we have for any open set $\Om\subset G$ the identity
\begin{equation}\label{automorphism}
\Tu_G(\varphi(\Om))=\frac {|\varphi(\Om)|}{|\Om|} \Tu_G(\Om)~.
\end{equation}
\end{corollary}

\begin{proof}
In our case $H=G$ and $\varphi$ is an automorphism. Clearly then
$K=\{0\}$ is the trivial group, $\mu_K=\delta_{0}$ is the trivial
measure, $K\cap\Om=\{0\}$, $\Tu_K(K\cap\Om)=1$,
$\mu_K(K\cap\Om)=1$ and $G/K\cong G$, $\mu_{G/K}\cong \mu_G$.
Thus we find $\nu=\mu_G\circ\varphi^{-1}$, and $C:= {d
\mu_H}/{d\nu}$ being constant, it can be computed on
$\Om^*:=\varphi(\Om)$ as $C=
{|\varphi^{-1}(\varphi(\Om))|}/{|\varphi(\Om)|}=
{|\Om|}/{|\varphi(\Om)|}$. Applying Proposition
\ref{prop:homomorphism} yields \eqref{automorphism} with $\le$
first. However, $\varphi^{-1}$ is also an automorphism, and that
implies the reverse inequality, too. Whence Corollary
\ref{cor:automorphinvariance} follows.
\end{proof}

The important special case when $G=\RR^d$ and $\varphi$ is any
linear mapping $A:\RR^d\to\RR^d$ was already noted in
\cite{arestov:hexagon}. There the computation of the constant $C$
is equivalent to the calculation of the volume element, ie. the
determinant, of the linear mapping $A$.

The next assertion was also observed in \cite{arestov:hexagon}
for $\RR^d$.

\begin{corollary}\label{cor:directproduct}
Let $G=G_1\times\dots\times G_n$ and $\Om_j\subset G_j \,\,
(j=1,\dots,n)$, $\Om=\Om_1\times\dots\times \Om_n$. Then we have
\begin{equation}\label{directproduct}
\Tu_G(\Om)=\Tu_{G_1}(\Om_1)\cdots\Tu_{G_n}(\Om_n) ~.
\end{equation}
\end{corollary}

\begin{proof}
The $\le$ direction easily follows from iteration of Proposition
\ref{prop:homomorphism}. On the other hand take any positive
definite functions $f_j$ on $G_j$ with $\supp f_j\subset\subset
\Om_j$ for $(j=1,\dots,n)$. It is easy to see that then the
product $f:=f_1\cdots f_n$ is a positive definite function on $G$,
with $\supp f\subset\subset \Om$, hence also the $\ge$ part of
\eqref{directproduct} follows.
\end{proof}

\mysubsection{Tur\'an constants on quotient groups}
\label{sec:quotient-groups}

\begin{corollary}\label{cor:quotientgroup}
Let G be a LCA group, K a closed subgroup of $G$, and suppose
that the Haar measures $\mu_K$ and $\mu_{G/K}$ of $G$ and $G/K$,
respectively, are normalized (as always) so that
$d\mu_G=d\mu_Kd\mu_{G/K}$. Let $\Om$ be any open set in $G$ and
$\Theta$ be its projection on $G/K$, ie.
$\Theta:=\{g+K\,:\,g\in\Om\}$. Then we have
\begin{equation}\label{quotient}
{\mathcal T}_G(\Omega) \le {\mathcal T}_{G/K}(\Theta){\mathcal
T}_K(\Omega\cap K)~~.
\end{equation}
In particular, if $\,\Om\cap K = \{0\}$, then ${\mathcal
T}_G(\Omega) \le {\mathcal T}_{G/K}(\Theta)$.
\end{corollary}
\begin{proof}
Consider $H:=G/K$ and the natural projection $\pi:G\to G/K$. It is
a continuous group homomorphism and thus Proposition
\ref{prop:homomorphism} can be applied with $\varphi:=\pi$. In
this case $\Theta=\pi(\Om)$ comprises the class of cosets $K+g$ so
that $K+g\cap\Om\ne \emptyset$, and the arising measure $\nu$ is
identical to $\mu_{G/K}$. Hence $C=1$ and we are led to
\eqref{quotient}. The special case is obvious.
\end{proof}

\mysubsection{Restrictions to subgroups and the Tur\'an
constants} \label{sec:sub-groups}

\noindent We show now that there is some sort of monotonicity in
the first argument of ${\cal T}_G(\Omega)$ as well.

\begin{corollary}\label{cor:subgroup}
Let G be a compact abelian group, and K a closed subgroup of $G$.
Let the Haar measures $\mu_K$ and $\mu_{G}$ be normalized
arbitrarily, and let $\Om$ be any open set in $G$. Then we have
\begin{equation}\label{generalsubgroup}
\Tu_G(\Omega) \le \frac{|G|}{\Abs{K}}\Tu_K(\Omega\cap K)~~.
\end{equation}
Here $|G|=\mu_G(G)$ and $|K|=\mu_K(K)$.
\end{corollary}

\begin{proof}
With $\mu_G$ and $\mu_K$ already given, we can define the Haar
measure $\mu_{G/K}$ so that the condition $\mu_G=\mu_K \mu_{G/K}$
still holds. Let $\varphi:=\pi$ and $H:=G/K$ as in the previous
Corollary. Since we always have $\Theta\subset G/K$, and thus
${\mathcal T}_{G/K}(\Theta)\le {\mathcal T}_{G/K}(G/K) = \mu_{G/K}(G/K)$, an
application of Corollary \ref{cor:quotientgroup} yields
$\Tu_G(\Om)\le \mu_{G/K}(G/K) \Tu_K(\Om\cap K)$. It remains to
see that for a {\it compact} group $G$ and (closed, hence
compact) subgroup $K$ also the quotient is compact, and according
to our choice of normalization we have
$\mu_{G/K}(G/K)=\mu_G(G)/\mu_K(K)$. The assertion follows.
\end{proof}

\begin{example}
Let us remark here that Lemma 1 of \cite{gorbachev:turan} can be
proved via Corollary \ref{cor:subgroup} by taking $G=\TT$,
$\Omega$ to be the interval $(-p/q,p/q) \subseteq \TT$ (for some
co-prime integers $p$ and $q$ with $p/q\le 1/2$) and $K$ to be
the (finite) subgroup of $\TT$ generated by $1/q$. The results in
\cite{gorbachev:turan} first show that the Tur\'an problem in
this case can be reduced to a finite problem of linear
programming (this is obviously the case for any Tur\'an problem
on a finite group) and Corollary \ref{cor:subgroup} shows half
the reduction. The reverse inequality is also true in this
particular case (this can be shown by ``convolving'' a positive
definite function on the subgroup with a Fej\'er kernel of
half-base $1/q$) but it cannot be expected to hold in general.
\end{example}


\mysection{Upper bound from packing}
\label{sec:upper-bound-from-packing}

\noindent
Here we show in \S\ref{sec:packing-proofs} the three main results
which give us upper bounds for the Tur\'an constant using ``packing''.
In the remaining part of this section we show several examples and
applications of these, in various groups.

\mysubsection{Proof of the main bounds from ``packing''}
\label{sec:packing-proofs}

\noindent
{\em Proof of Theorem \ref{th:diff}.}
Define $F:G\to\CC$ by
$$
F(x) = \sum_{\lambda,\mu\in\Lambda}f(x+\lambda-\mu).
$$
In other words $F = f*\delta_\Lambda*\delta_{-\Lambda}$, where
$\delta_A$ denotes the finite measure on $G$ that assigns a unit
mass to each point of the finite set $A$. It follows that $\ft{F}
= \ft{f}\Abs{\ft{\delta_\Lambda}}^2 \ge 0$ so that $F$ is
continuous and positive definite. Moreover, we also have
\begin{equation}
\supp F \subseteq \supp f + (\Lambda-\Lambda)
  \subseteq \Omega+(\Lambda-\Lambda)
\end{equation}
and
\begin{equation}\label{Ff-0}
F(0) = \Abs{\Lambda} f(0),
\end{equation}
since $\Omega \cap (\Lambda-\Lambda) \subseteq \Set{0}$.
Finally
\begin{equation}\label{Ff-1}
\int_G F = \Abs{\Lambda}^2 \int_G f.
\end{equation}
Applying the trivial upper bound $\int_G F \le F(0) \Abs{\Omega+(\Lambda-\Lambda)}$
to the positive definite function $F$ and using \eref{Ff-0} and \eref{Ff-1} we get
\begin{equation}\label{bound}
\int_G f \le {\Abs{\Omega+(\Lambda-\Lambda)} \over \Abs{\Lambda}} f(0).
\end{equation}
Estimating trivially $\Abs{\Omega+(\Lambda-\Lambda)}$ from above by $\Abs{G}$
we obtain the required ${\cal T}_G(\Omega) \le \Abs{G} / \Abs{\Lambda}$.
\Qed

\begin{corollary}\label{cor:packing}
Let $G$ be a compact abelian group and suppose $\Omega,H,\Lambda
\subseteq G$, $H+\Lambda\le G$ is a packing at level $1$,
that $\Omega \subseteq H-H$ and that $f\in\FF(\Om)$. Then
\eqref{bound-from-diff} holds.

In particular, if $H+\Lambda=G$ is
a tiling, we have
\begin{equation}\label{bound-for-tiles}
{\cal T}_G(\Omega) \le \Abs{H}.
\end{equation}
\end{corollary}

\begin{proof}
Since $H+\Lambda \le G$ it follows that $(H-H) \cap
(\Lambda-\Lambda) = \Set{0}$. Since $\Omega \subseteq H-H$ by
assumption it follows that $\Omega$ and $\Lambda-\Lambda$ have at
most $0$ in common. Theorem \ref{th:diff} therefore applies and
gives the result. If $H+\Lambda=G$ then
$\Abs{G}/\Abs{\Lambda}=\Abs{H}$ and this proves
\eref{bound-for-tiles}.
\end{proof}

\begin{proof}[Proof of Theorem \ref{th:diff-infinite}.]
Let $\epsilon>0$ and choose $R>0$ and $x\in G$ such that
$$
\Abs{\Lambda \cap Q_R(x)} \ge (\rho-\epsilon) \Abs{Q_R(x)} 
\ge (\rho-\epsilon)(R-1)^d,
$$
where $Q_R(x)$ is the cube of side $R$ and center at $x$.
Assume also that $\supp f \subseteq Q_r(0)$.

Let $\Lambda' = \Lambda \cap Q_R(x)$ and
construct the function $F$ as in the proof of Theorem \ref{th:diff},
with $\Lambda'$ in place of $\Lambda$.
We now have that
$$
\supp F \subseteq \supp f  + (\Lambda'-\Lambda') \subseteq Q_{2R+r}(0).
$$
This time we do not apply the trivial upper estimate to $F$ as we
did in Theorem \ref{th:diff} (then, we had no detailed
information on the support). Instead we use that for $L\in 2\NN$
\begin{equation}\label{non-trivial-bound}
{\cal T}_G(Q_L(0)) \le (L/2+1)^d~.
\end{equation}
The validity of ${\cal T}_{\RR^d}(Q_L(0)) \le 2^{-d}L^d ~ (\forall
L>0)$ and hence \eref{non-trivial-bound} in the case of $G=\RR^d$
has been proved, for example, in
\cite{arestov:hexagon,arestov:tiles,kolountzakis:turan}. For $G =
\ZZ^d$ we give a proof here.

Notice first that for any finite $\Omega \subseteq \ZZ^d$ and any
large enough positive integer $M$ we have
\begin{equation}\label{finite-sets-in-the-grid}
{\cal T}_{\ZZ^d}(\Omega) \le {\cal T}_{\ZZ_M^d}(\Omega).
\end{equation}
Indeed, if $M$ is large enough (e.g. $M>{\rm diam} (\Om)/2 $) then
the closed subgroup $K:=M\ZZ^d$ only intersects $\Om$ in $0$,
while the factor group $\ZZ_M^d$ will have an injective image
$\Theta$ of $\Om$: hence Corollary \ref{cor:quotientgroup}
applies.

If $\Omega=Q_L^d(0)=\Set{-L/2,\ldots,L/2}^d$ define $H$ to be the
set $\Set{0,\ldots,L/2}^d$ such that $\Omega=H-H$. Take now $M =
10(L/2+1)$, for example, so that (a) $H$ tiles $\ZZ_M^d$ by
translation, and, (b) $M$ is large enough to have all elements of
$\Omega$ distinct mod $\ZZ_M^d$. Using Corollary
\ref{cor:packing} we obtain \eref{non-trivial-bound} from
\eref{bound-for-tiles} in the group $\ZZ_M^d$,   and hence also
in $\ZZ^d$ because of \eref{finite-sets-in-the-grid}.

Hence taking $L:=L(R,r)$ in \eqref{non-trivial-bound} as the
least even integer not less than $2R+r$, we obtain both for
$G=\RR^d$ and $G=\ZZ^d$ the estimate $\int_G F \le \Tu_G(Q_L(0))
F(0) \le (R+r/2+2)^d F(0)$. Comparing this with \eref{Ff-0} and
\eref{Ff-1} (with $\Lambda'$ in place of $\Lambda$) we are led to
$$
\int_G f \le f(0) \frac{(R+r/2+2)^d}{\Abs{\Lambda'}}
 \le \frac{(R+r/2+2)^d}{(\rho-\epsilon)R^d}~.
$$
Since $\epsilon>0$ can be taken arbitrarily small and $R$
arbitrarily large, we get $\int_G f \le {1 \over \rho} f(0)$.
\end{proof}

\mysubsection{Sharpness}
\label{sec:sharpness}

\noindent The bound \eref{bound-from-diff} can be sharp. Take,
for example, $\Omega$ to be a subgroup of $G$ of finite index and
$H=\Omega$. Take also $\Lambda$ to a complete set of coset
representatives of $G/\Omega$, so that $\Abs{\Lambda}<\infty$.
Then $H+\Lambda=G$ and Corollary \ref{cor:packing} applies and
gives \beql{tmp-2} \sum_{x\in G} f(x) \le \Abs{\Omega}f(0) \eeq
for every positive definite function $f:G\to\CC$ supported in
$\Omega$, which is also the trivial bound. Taking
$f=\chi_\Omega$, which is positive definite because $\Omega$ is a
group, gives equality in \eref{tmp-2}.

More generally (and as in the next example) the inequality \eref{bound-from-diff}
is sharp whenever $H+\Lambda=G$ and $\Omega=H-H$.
In such a case the function $f = \chi_H*\chi_{-H}$ achieves equality in \eref{bound-from-diff}.

\mysubsection{Examples}

\begin{example}\label{ex:1}
Take $G=\ZZ_8=\Set{0,1,\ldots,7}$,
$H=\Set{0,1,4,5}$, $\Omega=H-H=\Set{0,1,3,4,5,7}$ and $\Lambda =
\Set{0,2}$, so that $\Lambda-\Lambda=\Set{0,2,6}$ and
$H+\Lambda=G$. It follows that
$$
\sum_{x\in G}f(x) \le 4 f(0)
$$
for any positive definite function on $\ZZ_8$ which vanishes on
$\pm2$, instead of the trivial $\sum_{x\in G}f(x) $ $\le 6 f(0)$.
The equality can be achieved by the function $f =
\chi_H*\chi_{-H}$.
\end{example}

\begin{example}\label{ex:2}
Let $G:=\ZZ$ and $\Om:=\Om_N:=\{-N,-1,0,1,N\}$; then the trivial
estimate is $A(N):=\Tu_{\ZZ}(\Om_N)\le 5$. Let $f\in\FF(\Om)$ be
a positive definite and real valued function: then $f(k)=f(-k)$,
that is, $f$ is even. The dual group is $\TT$, and positive
definiteness of $f$ means $p(x):=1+2f(1)\cos x + 2f(N) \cos Nx
\ge 0$ (as $f(0)=1$ by normalization). In the Tur\'an problem we
are to maximize $\int_\ZZ f =1+2f(1)+2f(N)=p(0)$; we have
$A(N)=\max p(0)$.

To find $A(N)$ in case when $N=2n+1$ is odd we may look at the
value $p(\pi)=1-2f(1)-2f(2n+1)\ge 0$ to see that $p(0)= 2 -
p(\pi) \le 2$. Clearly, any function with $f(1)+f(2n+1)=1/2$
achieves this bound while $p\ge 0$ if additionally we require
$0\le f(1),f(2n+1)$. Hence $A(2n+1)=2$.

If $N=2n$ is even, the solution is less simple. We claim that
$A(N)=1+1/\cos\frac{\pi}{2n+1}=:C(N)$, say, and the extremal
function is
$$
p_0(x):=1 + \frac {2n}{(2n+1)\cos \frac{\pi}{2n+1}} \cos x + \frac
{1}{(2n+1)\cos \frac{\pi}{2n+1}} \cos 2nx \,.
$$
Clearly $p_0(0)=C(N)$, and standard calculus proves nonnegativity
of $p_0$, hence it is an admissible trigonometric polynomial and
$A(N)\ge C(N)$.

To show its extremality we consider a general $p(x)=1+a\cos
x+b\cos 2nx$ (where $a:=2 f(1)$, $b:=2f(N)$) at the point
$z_0:=\pi+\pi/(2n+1)$, which yields $0\le
p(z_0)=1-a\cos\frac{\pi}{2n+1}-b\cos\frac{\pi}{2n+1}$. Thus
$p(0)=1+a+b = 1+(1-p(z_0))/\cos\frac{\pi}{2n+1} \le C(N)$, and
the calculation is concluded.

Now let us consider the estimates obtainable from the use of
Theorem \ref{th:diff-infinite}. In case $N$ is odd, taking
$\Lambda:=2\ZZ$ is optimal. Indeed, since $\Lambda$ is a
subgroup, $\Lambda-\Lambda=\Lambda$, and it does not intersect
$\Om_N$ (apart from $0$), hence an application of Theorem
\ref{th:diff-infinite} gives the right value $A(N)\le 1/{\rm
dens}(\Lambda)=2$. Hence in this case Theorem
\ref{th:diff-infinite} is sharp.

Let us see that it is {\it not} in the case when $N=2n$ is
even. To this, first we have to find the best upper density, that
is,
$$
L(N):=\sup_{\Om_N \cap (\Lambda-\Lambda) =\{0\}} \overline{{\rm
dens}}(\Lambda)\, .
$$

Let us consider the set
$\Lambda^*:=\{0,2,\dots,2n-2\}\cup\{2n+1,2n+3,\dots,4n-1\}+(4n+2)\ZZ$,
which contains $2n$ elements in each interval
$\big[k(4n+2),(k+1)(4n+2)\big)$ of $4n+2$ numbers and hence has
density $n/(2n+1)$. A direct calculation shows that $\Om_N \cap
(\Lambda^*-\Lambda^*) =\{0\}$, hence $L(N)\ge n/(2n+1)$. On the
other hand we assert that for no $\Lambda$ satisfying $\Om_N \cap
(\Lambda-\Lambda) =\{0\}$ can any interval $I=[k,k+2n]$ of $2n+1$
consecutive numbers contain more than $n$ elements of $\Lambda$.
Indeed, no pair of neighboring numbers belong to $\Lambda$,
because $1\in\Om_N$, and (at least) $n+1$ non-neighboring numbers
can be placed into $I$ only if all $m\in I$ with the same parity
as $k$ is contained. However, then both $k$ and $k+2n$ is
contained, having difference $2n\in\Om_N$, a contradiction. Hence
for a $\Lambda$ satisfying our condition, the upper density can
not exceed $n/(2n+1)$, which proves $L(N)=n/(2n+1)$.

Now we can compare the best estimate $\Tu_\ZZ(\Om_N)\le
1/L(N)=2+1/n$ arising from Theorem \ref{th:diff-infinite} to the
exact value $2+1/\cos \frac {\pi}{2n+1}$ found above. It shows
that application of Theorem \ref{th:diff-infinite} -- although
much better than the trivial estimate, but still -- is not
optimal in this case. This example highlights also the fact that
number theoretical, intrinsic structural properties -- like e.g.
$N$ being even or odd -- essentially influence the values of the
Tur\'an constants and sharpness of the estimates we have.
\end{example}

\begin{example}\label{ex:3}
Another example of a nice set with
nontrivial, but not sharp estimate arising from Theorem
\ref{th:diff-infinite} is the unit disk $D$ in $\RR^2$ (with Lebesgue measure).
The area of $D$ is $\pi$ and the right value of the Tur\'an constant,
first computed by Gorbachev \cite{gorbachev:ball}, is
$|D|/2^d=\pi/4$ in this case. Now $D$ is the difference set of
$H:=D/2$, and the best density we can have is, in fact, the {\em
sphere packing constant} of $\RR^2$. It is well-known \cite{PA}
that the best packing is the regular hexagon lattice packing, hence
$L(D)=2/\sqrt{3}$ and the arising estimate is $\sqrt{3}/2$. In
comparison, note that the estimate of \S \ref{sec:dispersed}
gives $|D|/2=\pi/2$, while the estimate of Theorem
\ref{th:spectral-infinite} from the spectral approach does not
apply, since the ball is {\it not} spectral. The above values
compare as
$\pi/4=0.785\dots<\sqrt{3}/2=0.866\dots<\pi/2=1.57\dots$.
\end{example}

\begin{example}\label{ex:3.5}
We see that for a general $\Om\subset H-H$ or even $\Om=H-H$ the
``best translational set'', (i.e. the maximal number of elements or
the highest possible upper density), does not always achieve an
exact bound of $\Tu_G(\Om)$. In this respect it is worth
mentioning that, on the other hand, results of Herz
\cite{herz:example1}, \cite{herz:example2} show that each
subgroup $\Lambda$ of $G$ provides the theoretically best
possible, sharp estimate for {\it some} open set $\Om$. E.g. if
$G$ is compact, and $\Lambda$ is a finite subgroup having $n$ elements,
there exists a Borel set $H$ with the properties $|H|=1/n$,
$\Om:=H-H$ is open, and $\Om\cap\Lambda=\{0\}$. See also
\cite[7.4.1]{rudin:groups}. Clearly for this $\Om$ and $H$
we have that $H+\Lambda=G$ is a tiling, and $\Tu_G(\Om)=1/n$, achieved by
$\chi_H*\chi_{-H}$.
\end{example}


\begin{example}\label{ex:4}
The size of the Tur\'an constant of a set $\Omega$ may be
extremely small. Take for example in the group $G=\ZZ_{2n}$ the
set $\Omega=\Set{0} \cup K^c$, where $K$ is the subgroup generated
by $2$. Let then $\Lambda=K$ and apply Theorem \ref{th:diff}. It
follows that ${\cal T}_G(\Omega) \le 2$ while $\Abs{\Omega}=n+1$.

The same way we have $\Tu_{\ZZ} (\Om) \le 2$ for any subset
$\Om\subset (\{0\}\cup (2\ZZ+1))$ in view of Theorem
\ref{th:diff-infinite} and considering the set $\Lambda:=2\ZZ$.
(This covers the $N$ odd case of Example \ref{ex:2}, too.)

The generality of this example should be obvious.
\end{example}

\mysubsection{The Tur\'an constant of difference sets of tiles in $\RR^d$ or $\ZZ^d$.}

Here we show how to generalize the results in
\cite{arestov:tiles} (see also \cite{kolountzakis:turan}).
In \cite{arestov:tiles} the Tur\'an
constant of convex polytopes which tile $\RR^d$ by lattice
translation was determined.

Actually being a polytope and {\em lattice} translation need not
be assumed as it is a fact (see e.g. the references in \cite{kolountzakis:turan})
that any convex body that tiles space by translation is a
polytope and can also tile by lattice translation.

From Theorem \ref{th:diff-infinite} it follows that if $H$ is any measurable set
of finite measure that tiles $\RR^d$ or $\ZZ^d$ by translation with $\Lambda$
then the Tur\'an constant
of $H-H$ is equal to $1/\dens\Lambda = \Abs{H}$.

Whenever $\Omega$ is a convex body in $\RR^d$ one can take
$H={1\over 2}\Omega$, so Theorem \ref{th:diff-infinite} is indeed
a generalization of the result in \cite{arestov:tiles}.

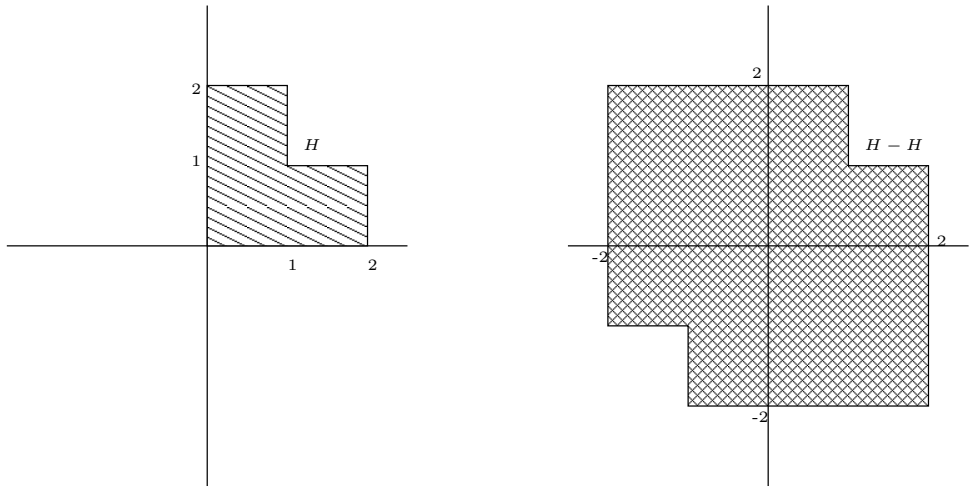
\begin{figure}
\begin{center}\input{example.pstex_t}\end{center}
\caption{The Tur\'an constant of $H-H$ (right) is equal to the area of $H$}
\label{fig:domains}
\end{figure}

However, Theorem \ref{th:diff-infinite} can determine the Tur\'an constant
of many more sets than those dealt with in \cite{arestov:tiles},
such as the one in Figure \ref{fig:domains}.
The subset of the plane $H$ shown on the left tiles the plane by
translation hence its difference set shown on the right has
Tur\'an constant equal to $\Abs{H}$.

\begin{example}\label{ex:5}
Let $H\subset \ZZ^2$ be the three-element set
$\{(0,0),\,(0,1),\,(1,0)\}$ and $\Om$ be the difference set
$H-H=\{(-1,0),\,(-1,1),\,(0,-1),\,(0,0),\,(0,1),\,(1,-1),\,(1,0)\}$.
Then $|\Om|=7$, but $H$ tiles $\ZZ^2$, hence Theorem
\ref{th:diff-infinite} applies and yields $|H|=3$.
Observe that the set $\Lambda:=\ZZ (1,1) + \ZZ (2,-1)$ provides a
translational set. Indeed, any points $(n+2m,n-m)$ of $\Lambda$,
and thus also of $\Lambda-\Lambda$, has the property that the
first coordinate is congruent to the second $\mod 3$, hence
$\Om\cap\Lambda-\Lambda=\{(0,0)\}$. On the other hand all points
of $\ZZ^2$ with the above congruence property belong to
$\Lambda$, i.e. $\Lambda$ is a subgroup of index $3$. It follows
that the density of $\Lambda$ is $1/3$, and Theorem
\ref{th:diff-infinite} gives the assertion.
\end{example}

\mysubsection{The Tur\'an constant of dispersed
sets}\label{sec:dispersed}

\noindent As an application of Theorem \ref{th:diff-infinite} we
show that, in $\RR$, the Tur\' an constant of a set of given
length is the largest if the set is an interval. The construction
extends to $\ZZ$, and even to $\RR^d$ and $\ZZ^d$ giving a
generally valid improvement of the trivial bound by about a factor
of $2$.

\begin{theorem}\label{th:d-dim-interval}
Let $\Omega\subseteq \RR^d$ be an open set of finite measure $m$.
Then we have
\begin{equation}\label{d-dim-turan-constant-bound}
{\cal T}_{\RR^d}(\Omega) \le \frac m2 ~.
\end{equation}
Let $\Omega\subseteq\ZZ^d$ be a set of size $m$ containing the
origin and denote by $m^{+}$ the number of lattice points in the
"nonnegative half of $\Om$", i.e. in
$\Omega\cap\big([0,\infty)\times\ZZ^{d-1}\big)$. Then we have
\begin{equation}\label{Zd-turan-constant-bound}
{\cal T}_{\ZZ^d}(\Omega) \le m^{+}~.
\end{equation}
\end{theorem}
\begin{proof}
Let us denote $P:=[0,\infty)\times\RR^{d-1}$ or
$[0,\infty)\times\ZZ^{d-1}$, respectively, and put $\Om^+:=\Om\cap
P$. Note that in $\RR^d$ we simply have $m^{+}:=|\Om^+|=m/2$. It
is easy to see that Theorem \ref{th:equivalences} (on the
equivalent formulations of the Tur\'an constant), allows us to
assume that $\Om$ is bounded: so let $\Om\subset B(0,r)$ with
some fixed ball of radius $r$. Take a large parameter
$L_0>\max\{2,r\}$, define $L_k=L^{2^k}=L_{k-1}^2$ ($\forall
k\in\NN$), say, and put
\begin{equation}\label{def-of-Q}
Q_k:=Q_{L_k}((L_k,0,\dots,0))=[0,2L_k]\times[-L_k,L_k]^{d-1} \quad
(k\in \NN),\quad Q_0:=\emptyset ~.
\end{equation}
Note that $|Q_k|=(2L_k)^d$ in $\RR^d$ and $(2L_k+1)^d$ in $\ZZ^d$.
Define
\begin{equation}\label{def-of-S}
S_k:=Q_k\setminus (Q_{k-1}+\Om) \quad (k\in \NN) ~.
\end{equation}
Obviously, $S_k$ are closed sets of measure
\begin{equation}\label{measure-of-S}
|S_k|\ge |Q_k|-|Q_{k-1}+\Om| \ge (2L_k)^d-((2L_{k-1}+1)+2r)^d \ge
2^d L_k^d\left(1-\big(\frac{2+r}{L_{k-1}}\big)^d\right) \quad
(k\in \NN),
\end{equation}
satisfying $(S_k-S_n)\cap \Om = \emptyset$ for $k\ne n$. We aim at
constructing the discrete set
\begin{equation}\label{def-of-Lambda}
\Lambda:=\bigcup_{k=1}^{\infty} \Lambda_k,\qquad \Lambda_k\subset
S_k \quad (k\in \NN) ~
\end{equation}
with as many as possible elements but satisfying
$(\Lambda_k-\Lambda_k)\cap \Om = \{0\}$. Note that if the latter
condition is satisfied, then we will also have
$(\Lambda-\Lambda)\cap \Om = \{0\}$ in view of the respective
property of $S_k\supset \Lambda_k$. So now we define the elements
of $\Lambda_k$ inductively by a ``greedy algorithm'' as follows.
Let $\lambda_0^{(k)}$ be any element of the nonempty set $S_k$
with first coordinate $0$. Such an element clearly exists. Then
for $n\ge 1$ take any
\begin{align}\label{def-of-lambda-kn}
\lambda_n^{(k)}&:=(x_{1,n},\dots,x_{d,n})\in \big( S_k\setminus
\bigcup_{j=1}^{n-1} (\lambda_j^{(k)}+\Om^+)\big)
\notag \\ &{\rm with} \\
x_{1,n}&=\min\bigg\{x_1~:~ \exists ~ x=(x_1,\dots,x_d)\in \big(
S_k\setminus \bigcup_{j=1}^{n-1} (\lambda_j^{(k)}+\Om^+)\big)
\bigg\}~.\notag
\end{align}
Defining new elements $\lambda_n^{(k)}$ of $\Lambda_k$ terminates
in a finite number of steps, but not before $\cup_{j=1}^{n-1}
(\lambda_j^{(k)}+\Om^+)$ covers $S_k$, so with $m^{+}:=|\Om^{+}|$
we must have
\begin{equation}\label{size-of-Lambda}
\# \Lambda_k \ge \frac{|S_k|}{|\Om^+|}\ge \frac{2^d
L_k^d\left(1-(\frac{2+r}{L_{k-1}})^d\right)}{m^{+}} \quad (k\in
\NN) ~.
\end{equation}
By construction, for any $n>j$
$\lambda_n^{(k)}-\lambda_j^{(k)}\in\Om$ is not possible, hence
$\Lambda-\Lambda\cap \Om =\{0\}$. Moreover, in view of
\eqref{size-of-Lambda} we have
\begin{equation}\label{density-of-Lambda}
\overline{\dens} \Lambda \ge \limsup_{k\to \infty} \frac{\#
\Lambda_k}{|Q_k|} \ge \limsup_{k\to \infty}
\frac{\left(1-(\frac{2+r}{L_{k-1}})^d\right)}{m^{+}} =\frac
1{m^{+}}~.
\end{equation}
Now an application of Theorem \ref{th:diff-infinite} with
$\Lambda$ concludes the proof.
\end{proof}

\begin{remark} For $d=1$ \eqref{d-dim-turan-constant-bound}
is sharp for intervals in $\RR$. It is plausible, but we do not
know if intervals are the only cases of equality.
\end{remark}

\begin{remark} As $\Omega$ is always symmetric, in $\ZZ$ we always have
$m^{+}=(m+1)/2$. The estimate \eqref{Zd-turan-constant-bound} can
also be sharp at least for $d=1$. Take e.g. $\Omega=\Om_0$ or
$\Om_1$ from Example \ref{ex:2}, or, more generally, take
$\Om:=[-N,N]$. Then $m=2N+1$, $m^{+}=(m+1)/2=N+1$, and the Fej\'er
kernel shows that this value can be achieved. Thus
$\Tu_{\ZZ}([-N,N])=N+1$, and intervals have maximal Tur\'an
constants once again. However, here the sets
$k[-N,N]:=\{kn~:~|n|\le N\}$ of similar size have equally large
Tur\'an constants, hence intervals are not the only extremal
examples in $\ZZ$.
\end{remark}

\begin{remark} It can be proved that the asymptotic uniform upper
density of all sets remain the same both in $\RR^d$ and in
$\ZZ^d$ if we define it replacing $Q_R$ by $RK$ with any other
convex body $K$. Thus in the above proof one can consider the
slightly modified basic sets $RQ(T)$, where $RQ(T)$ is the
$R$-dilated copy of the unit box rotated by the isometry $T\in
SO(d)$. If we choose $T$ to be "irrational" in the sense that no
lattice point (apart from the origin) moves to the hyperplane
$\{x_1=0\}$, then with these sets a similar argument leads the
same estimate but now with $m^{+}=\# \Om^{+}=(m+1)/2$. We leave
the details to the reader.
\end{remark}

\mysubsection{The Tur\'an constant of an interval missing two points}

\noindent Our next result shows the effect of forcing a positive
definite function to vanish at a neighborhood of one point in an
interval.
\begin{theorem}\label{th:interval-minus-point}
Suppose $0 < b < a \le 2b$ and let
$$
\Omega = (-a,-b) \cup (-b,b) \cup (b,a).
$$
Then ${\cal T}_\RR(\Omega) = {\cal T}_\RR(-b,b) = b$.
\end{theorem}

\begin{proof}
Simply take $\Lambda = b\ZZ$ and apply Theorem \ref{th:diff-infinite}
to obtain that ${\cal T}_\RR(\Omega) \le  b$.
The other direction is obvious by the monotonicity of ${\cal T}_G(\cdot)$.
\end{proof}

The condition $a<2b$ is necessary in Theorem
\ref{th:interval-minus-point}. Indeed, if $a>2b$ then, with
$c=\min\{b,(a-b)/2\}>b/2$ and $d:=(a+b)/2$ the function
$f:=\chi_{(0,c)}*\chi_{(-c,0)}*(\delta_{0}+\delta_{d})
*(\delta_{0}+\delta_{-d})$, whose graph consists of three
triangles centered at $0$ and $\pm d$ of width $2c$ and heights
$1$ (for the central triangle) and $1/2$ (for the other two) is
positive definite and supported in $\Omega$, yet has $f(0)=2c$
and $\int_\RR f = 4 c^2$. Hence $\Tu_\RR (\Om)\ge 2c>b$.

\mysection{Upper bound from spectral sets}

\mysubsection{Proof of the bound from spectral sets}
\label{sec:upper-bound-from-spectral-sets}

\noindent
{\em Proof of Theorem \ref{th:spectral}.}
Since $T$ is a spectrum of $H$ we have (see \S\ref{sec:spectral})
\begin{eqnarray*}
\supp\ft{\chi_T} &\subseteq&  \Set{0} \cup (H-H)^c\\
  &\subseteq&  \Set{0} \cup \Omega^c\\
  &\subseteq&  \Set{0} \cup \Set{f=0}.
\end{eqnarray*}
Hence $\ft{f}+T=c\ft{G}$ is a tiling and $c=\Abs{T}f(0)$,
as $\int_{\ft{G}}\ft{f} = \Abs{\ft{G}}f(0)$.

Since $\ft{f}\ge 0$ in $\ft{G}$ it follows that $\ft{f}(0) \le c$ or
$$
\sum_{x\in G}f(x) \le \Abs{T}f(0) = \Abs{H}f(0).
$$
\Qed

\mysubsection{Comparison of Theorems \ref{th:diff},
\ref{th:diff-infinite} and \ref{th:spectral},
\ref{th:spectral-infinite}} \label{sec:comparison}

\noindent First we give an example when Theorem \ref{th:spectral}
gives a better bound than any possible application of Theorem
\ref{th:diff}. Let $G = \ZZ_2^{12}$ and $H = \Set{e_1, e_2,
\ldots, e_{12}}$, where $e_i$ is the vector in $G$ with all zeros
except at the $i$-th position where we have $1$. The set $H$ was
recently shown by Tao \cite{tao:fuglede-down} to have a spectrum,
and it is clear that $H$ cannot tile $G$ since $\Abs{H}=12$ does
not divide $\Abs{G}=2^{12}$.

Let $\Omega=H-H$. This means that $\Omega$ consists of the
all-zero vector plus all vectors in $G$ with precisely two $1$'s,
hence $\Abs{\Omega} = {12 \choose 2}+1 = 67$.

By Theorem \ref{th:spectral} we have that if $f:G\to\CC$ is a positive definite function
supported on $\Omega$ then
$$
\sum_{x\in G}f(x) \le 12 f(0).
$$

Suppose now that Theorem \ref{th:diff} applies with some $\Lambda\subseteq G$,
such that $\Omega \cap (\Lambda-\Lambda) = \Set{0}$.
Since $\Omega=H-H$ this implies that $H+\Lambda\le G$ is a packing at level $1$,
hence $\Abs{\Lambda} \le {1\over 12}\Abs{G}$.
In fact $\Abs{\Lambda} < {1\over 12}\Abs{G}$ as $\Abs{\Lambda}$ is an integer
but ${1\over 12}\Abs{G}$ is not.
Clearly then \eref{bound-from-diff} is inferior than $\sum_{x\in G}f(x) \le 12 f(0)$
given by Theorem \ref{th:spectral}.

Tao \cite{tao:fuglede-down} also shows how to construct a domain
(in fact, a finite union of unit cubes) in $\RR^d$, $d\ge 5$,
which is spectral but not a translational tile. Suppose $H$ is
such a domain. Theorem \ref{th:spectral-infinite} shows that
${\cal T}_{\RR^d}(H-H) \le \Abs{H}$. We claim that Theorem
\ref{th:diff-infinite} gives a worse upper bound for the set
$\Omega=H-H$. Indeed, suppose that $\Lambda\subseteq\RR^d$ is a
set for which
\begin{eqnarray*}\label{omega-lambda}
\Omega \cap (\Lambda-\Lambda) = \Set{0},
\end{eqnarray*}
as required by Theorem \ref{th:diff-infinite}, and that $\rho$ is
the upper density of $\Lambda$. Condition \eref{omega-lambda}
means that $H + \Lambda$ is a packing, hence $\Abs{H}\dens\Lambda
\le 1$. The fact that $H$ is not a tile implies (this requires a
proof, an easy diagonal argument) that the inequality above is
strict, so that $1/\rho > \Abs{H}$, which shows that any
application of Theorem \ref{th:diff-infinite} gives a worse result
than Theorem \ref{th:spectral-infinite} for $H-H$.

\noindent
\ \\
{\bf Bibliography}
\\

\noindent
{\sc\small
Department of Mathematics, University of Crete, Knossos Ave., \\
714 09 Iraklio, Greece.\\
E-mail: {\tt kolount@member.ams.org}\\
\ \\
and\\
\ \\
{\sc\small
Alfr\' ed R\' enyi Institute of Mathematics, Hungarian Academy of Sciences, \\
1364 Budapest, Hungary}\\
E-mail: {\tt revesz@renyi.hu}
}

\end{document}

%% file: example.pstex_t
\begin{picture}(0,0)%
\includegraphics{example.pstex}%
\end{picture}%
\setlength{\unitlength}{4144sp}%
\begingroup\makeatletter\ifx\SetFigFont\undefined%
\gdef\SetFigFont#1#2#3#4#5{%
  \reset@font\fontsize{#1}{#2pt}%
  \fontfamily{#3}\fontseries{#4}\fontshape{#5}%
  \selectfont}%
\fi\endgroup%
\begin{picture}(5778,2901)(6406,-3379)
\put(8096,-2072){\makebox(0,0)[lb]{\smash{\SetFigFont{6}{7.2}{\familydefault}{\mddefault}{\updefault}{\color[rgb]{0,0,0}1}%
}}}
\put(8576,-2072){\makebox(0,0)[lb]{\smash{\SetFigFont{6}{7.2}{\familydefault}{\mddefault}{\updefault}{\color[rgb]{0,0,0}2}%
}}}
\put(7521,-1449){\makebox(0,0)[lb]{\smash{\SetFigFont{6}{7.2}{\familydefault}{\mddefault}{\updefault}{\color[rgb]{0,0,0}1}%
}}}
\put(7521,-1017){\makebox(0,0)[lb]{\smash{\SetFigFont{6}{7.2}{\familydefault}{\mddefault}{\updefault}{\color[rgb]{0,0,0}2}%
}}}
\put(10877,-2983){\makebox(0,0)[lb]{\smash{\SetFigFont{6}{7.2}{\familydefault}{\mddefault}{\updefault}{\color[rgb]{0,0,0}-2}%
}}}
\put(11980,-1929){\makebox(0,0)[lb]{\smash{\SetFigFont{6}{7.2}{\familydefault}{\mddefault}{\updefault}{\color[rgb]{0,0,0}2}%
}}}
\put(10877,-922){\makebox(0,0)[lb]{\smash{\SetFigFont{6}{7.2}{\familydefault}{\mddefault}{\updefault}{\color[rgb]{0,0,0}2}%
}}}
\put(9918,-2024){\makebox(0,0)[lb]{\smash{\SetFigFont{6}{7.2}{\familydefault}{\mddefault}{\updefault}{\color[rgb]{0,0,0}-2}%
}}}
\put(8192,-1353){\makebox(0,0)[lb]{\smash{\SetFigFont{6}{7.2}{\familydefault}{\mddefault}{\updefault}{\color[rgb]{0,0,0}$H$}%
}}}
\put(11549,-1353){\makebox(0,0)[lb]{\smash{\SetFigFont{6}{7.2}{\familydefault}{\mddefault}{\updefault}{\color[rgb]{0,0,0}$H-H$}%
}}}
\end{picture}